# Alice Ambrose on Logic, A Priori Concepts, and the Epistemology of Convention

Author: Juan J. Colomina-Alminana

Affiliation: Northeastern University, Oakland Campus (Mills College)

Address: 5000 MacArthur Blvd., 132 Vera Long, Oakland, CA, 94613 – USA

Email: j.colomina-alminana@northeastern.edu

Declaration: The author declares that there is no conflict of interest.

No data employed

This article adheres to all ethical standards

**Abstract:** This essay argues that Alice Ambrose's early published articles precede key elements of W.V. Quine's later critiques, first of "Truth by convention" (Quine 1936) and later of the analytic-synthetic distinction (Quine 1951). I demonstrate how Ambrose identifies in writing as early as in 1931: (1) the paradox of treating logical principles as mere conventions, (2) the infinite-regress in stipulative definitions, (3) the pre-conditional role of logic for any convention, and (4) the instability of the analytic-synthetic divide. Ambrose's insights, therefore, prefigure Margaret Macdonald's 1934 unpublished dissertation (Cf. Spinney 2025) and predate some of Quine's major contributions a few years before he made them popular.

**Keywords:** Early Analytic Philosophy, Logic, Apriority, Truth by Convention, Infinite Regress.

* * *

The received history of analytic philosophy often locates the critical turn away from conventionalism in W.V. Quine's landmark essay "Truth by Convention" (1936). There, Quine attacks the idea that logical truths are grounded on definitions or stipulations. Later, in "Two

Dogmas of Empiricism" (1951), Quine famously tries to dismantle the analytic-synthetic distinction altogether. Even though the narrative regarding when and how Quine challenges the notion of analyticity varies, it has only recently disputed its origin (Spinney 2025), and credit has also rightly been given to Margaret MacDonald in anticipating some of Quine's crucial arguments in her 1934 unpublished dissertation.

Unacknowledged, though, still are Alice Ambrose's early published critiques. Ambrose's philosophical writings from the 1930s, especially her 1931 criticism of C. I. Lewis's intensional logic system, her 1932 Ph.D. dissertation, partially published in 1933, and her 1935 defense of a moderate finitism, already contain deep skepticism about the epistemic grounding of a priori concepts and analytic truths, a methodological strategy that prefigures key elements of both Macdonald's and Quine's later positions.

This paper goes further. I argue that specifically in her 1931 published 17-page critical notice of C. I. Lewis's *Mind and the World Order* and her 1933 examination of the debate among intuitionists, formalists, and logisticians on the foundations of logic and mathematics, Ambrose had already laid bare the self-undermining nature of treating logic as a legislative fiat, the circularity inherent in stipulative definitions, and the necessity of logic as a precondition for any meaningful convention. In what follows then, I develop a detailed, one to one critical chronological comparison and overview of Ambrose's arguments with those of Macdonald and Quine, demonstrating that Ambrose anticipated both the critique of "truth by convention" found in MacDonald and Quine, and already problematizes the analytic-synthetic divide, which Quine will later fully attack.

One small note before moving forward. This paper argues that intellectual labor comes in association, and that it is very possible that clever minds think alike. After all, as Spinney (2025) suggests, many of these criticisms are rooted in the work of even earlier logicians. Therefore, this paper defends that, as the broad scope of their respective work shows, all three authors were arguing against conventionalism in logic while pursuing their individual and independent research agendas. Nevertheless, this article also demonstrates that as early as in 1931 Alice Ambrose had published arguments later discussed by both MacDonald and Quine.

Even though I will not engage with it fully, I want to highlight that a recent and valuable contribution to the recovery of Ambrose's place in early analytic philosophy is Connell's (2022) study, which examines Ambrose's later work in the philosophy of language and metaphilosophy, especially her development of themes drawn from Wittgenstein's middle period. While Connell's paper significantly advances our understanding of Ambrose's post-war contributions, it does not address the earlier logical and epistemological writings that are the focus of the present study.[1] My concern here is with Ambrose's 1931–1934 publications on logic, extensionalism, and the epistemology of convention—texts that predate the period Connell discusses and that, I argue, anticipate key elements of Macdonald's and Quine's later critiques.

The structure of the paper is as follows. Section 1 gives a brief presentation of both Quine's resistance to accept that logical truth is established by convention and Quine's criticism of the analytic-synthetic dichotomy. Section 2 offers a short overview of Ambrose's early publications. Presenting a one-to-one chronological comparison of the main arguments, Sections 3-4 closely analyze how Ambrose's early work contributes to the debate on truth by convention and anticipates both Macdonald's and Quine's critical insights. In addition, primarily by posing the infinite regress argument while defending the benefits of extensional logic over formal and intensional systems, Ambrose (1931 and 1933) also precede some of Quine's (1951) criticism of the analytic-synthetic dichotomy. Section 5 discusses some personal connections and affinities between Ambrose and MacDonald. Section 6 offers some concluding remarks.

---

[1] Two other important contributions are Chapman (2024) and Lonner (2024), which will be addressed later.

## 1. The Received Narrative

The standard historical narrative treats Quine's critique of conventionalism as unfolding in two distinct but often conflated stages. The first concerns the status of logical truths and the viability of grounding them in stipulation or linguistic convention. The second concerns the analytic–synthetic distinction and the broader epistemological architecture of empiricism. Although these two lines of argument are sometimes presented as continuous, they are in fact logically independent, and Quine himself develops them with different targets and different motivations. Clarifying this distinction is essential for understanding how Ambrose's early work anticipates elements of both.

In "Truth by Convention" (1936), Quine challenges the idea—common among logical positivists and present in C. I. Lewis—that logical truths can be rendered true solely by stipulation or definition.[2] His objections are primarily methodological and logical, not semantic. For two difficulties dominate. On one hand, one find the regress/circularity problem: Any attempt to derive logical truths from conventions presupposes the very logical principles needed to apply those conventions. As Quine puts it, deriving logic from convention requires logic itself, generating either circularity or an infinite regress. On the other hand, we have the supertask problem: To secure all logical truths by stipulation would require a complete specification of meaning postulates or translation rules—an impossible task, and one that again would presuppose logical inference. These objections target conventionalism in the foundations of logic and mathematics, not the analytic–synthetic distinction, for they concern the source of logical necessity but not the nature of meaning.

[2] There is substantial secondary literature, though, that emphasizes that Quine's 1936 position does not yet contain the full-blown rejection of analyticity found in 1951. Scholars such as Yemima Ben-Menahem (2005) and Greg Frost-Arnold (2013), for example argue that in "Truth by Convention" Quine still allows that one may add analytic stipulations to a linguistic framework, and that his later view—that one may also subtract or revise such stipulations—emerges only under the influence of Tarski's semantic work in the 1940s. Nothing in the present paper turns on denying this developmental trajectory, for my claim here is only that Ambrose anticipates several of the foundational objections Quine raises in 1936, not that she anticipates the later semantic holism of "Two Dogmas."

By contrast, “Two Dogmas of Empiricism” (1951) advances a different and broader thesis. There, Quine argues that the analytic–synthetic distinction itself lacks a non-circular foundation. His critique, thus, does not depend on the earlier regress argument. Instead, it rests on the instability of synonymy, the failure of reductionism, and the holistic character of confirmation. The conclusion is that analyticity cannot be defined without presupposing the very semantic notions it is meant to explain. This time he is emphasizing a semantic and epistemological critique but not a critique of conventionalism about logical axioms. Thus, while both essays undermine sharp boundaries between convention and empirical content, they do so via different argumentative routes. For the former attacks the foundational role of convention in logic whereas the latter attacks the semantic coherence of analyticity.

Understanding the distinction matters, for recognizing the independence of these two strands allows us to see more clearly what Ambrose anticipates. Her early writings (1931–1933) directly engage the foundational problem: the regress of stipulative definitions, the circularity of intensional systems, and the pre-conditional role of logic for any convention. These are precisely the themes Quine later develops in 1936. At the same time, Ambrose’s insistence on extensional grounding, her critique of intensional definitions, and her sensitivity to the instability of purely conceptual relations already gesture toward the semantic worries that Quine will later mobilize against analyticity in 1951. But these anticipations occur without collapsing the two projects, for Ambrose's arguments against conventionalism in logic stand on their own. Her challenges to intensional meaning and a priori concepts form a separate, though related, line of critique.[3]

## 2. Ambrose's Early Work.

Ambrose's early work is summarized in the following lemma: A very clever plea for a finitist extensional logic.[4] Alice Ambrose's (1931) core objection to an intensional logic is that, by defining every term solely by its relations to other concepts, you would end up either with viciously circular definitions or an infinite regress—and, worse, you have no way of reconnecting those purely conceptual relations to a world that, supposedly, we all share. An extensional logic breaks such a circle by fixing each predicate in the public realm of objects it picks out, anchoring meaning in observable referents and, thereby, guaranteeing the common ground needed for rigorous mathematics as well as genuine communication (Cf. Ambrose 1932).

In addition, Ambrose (1933) argues that unrestricted infinitary reasoning—employing "all" and "there exists" over infinite domains—commits us to the universal validity of the law of excluded middle, the solvability of every problem, and non-constructive reductio proofs even when no finite decision procedure exists, thus entangling us in circular definitions and other paradoxes inherent to infinity. By contrast, finitism confines quantification to finite, effectively calculable processes, restoring transparency in definitions, avoiding infinite regress, and ensuring every assertion carries an explicit finite verification or construction (Cf. Ambrose 1935a, 1935b, and 1938). Even though both sides are part of the same methodological and critical coin, and often one feeds the other, in this paper I argue that Ambrose's explicit defense of extensional logic (which was often the consensus at the time) carries important secondary caveats that anticipate criticism later proposed by Margaret MacDonald (1934) and W.V. Quine (1936 and 1951).[5]

[3] It should be noted that, as mentioned in the past note, the scholarly consensus has long recognized that Quine's arguments against analyticity did not emerge fully formed in 1936, nor are they attributed to him alone. Ben-Menahem (2005) and Frost-Arnold (2013), among others, document both the evolution of Quine's views and the influence of Tarski, Carnap, and other contemporaries. What I am arguing here is simply that there is a historiographical tendency to treat Quine's 1936 and 1951 critiques as the canonical turning points in the rejection of conventionalism and analyticity, but not the stronger claim that Quine was the sole originator of these lines of argument. For, as mentioned above, clever minds often think alike.

[4] Ambrose's early work can be divided in two separate groups. Ambrose 1931, 1932, 1933, 1934, and 1936a support extensional logic, such as the system developed in *Principia*, against the dangers of intension and the foundation of logic and mathematics on a priori concepts and stipulations. Ambrose 1935a, 1935b, 1936b, 1937, and 1938 directly address the necessity of finitism in logic. Even though both are part of the same endeavor, one clearly defined early on before her arrival in Cambridge in summer 1932, the second group of papers, as it has been argued by Connell (2022) and Loner (2024), seem at least influenced by some of Wittgenstein's moving thoughts during her stay in Cambridge (1932-1935). The first group though is clearly independent and, I argue, anticipates crucial arguments in contemporary philosophy.

In her 1931 paper "A Critical Discussion of *Mind and the World-Order*," Ambrose mounts a 17-page detailed criticism of C.I. Lewis's epistemological framework as explained in his book, which he defines in terms of Conceptualistic Pragmatism. While Lewis maintained that analytic truths are not falsifiable but could be pragmatically discarded, Ambrose questions the stability and epistemic legitimacy of such a stance. She writes: "Concepts may be discarded but cannot change… To be plausible, a priori categories and concepts... cannot 'float in from nowhere'" (Ambrose 1931: 367-368). Ambrose also challenges Lewis's narrow pragmatic conception of truth. She notes: "The pragmatic conception of truth is not exhausted by an account in which the whole test of truth is its achieving intelligible order" (Ambrose 1931: 377). Her resistance to ground logic on a priori concepts and her insistence on empirical verification and proof, therefore, anticipate Quine's later argument against truth by convention (Quine 1936) and the claim that the analytic-synthetic divide masks the revisability of all claims, including logical ones, under pressure from experience (Quine 1951).

[5] As mentioned, Ambrose's endeavor was already set before arriving in Cambridge. For example, some of the engagement with ordinary language, at least as contraposed to logical and mathematical language and propositions, is already shown in Ambrose (1933). She writes: "I have indicated now the bearing of a restricted right of definition on such propositions as Zermelo's axiom of inclusion and on the definition of the continuum. The crux of the issue is such terms as 'all' and 'there exists.' […] The formalist holds that… '[n]ot all' is thus made equivalent to 'there exists,' whether or not an infinite logical sum is entailed in the object's discovery" (1933: 602-603). This also anticipates what we now known as Quinean ontological commitment: The principle that a theory is committed only to the entities over which it quantifies.

Ambrose's 1932 unpublished dissertation investigates extensionalism in logic as developed in systems such as that in Whitehead and Russell's *Principia*. There she emphasizes that logical systems must be epistemically and semantically grounded, not merely formal. Reflecting on this, she writes later in life: "Investigation of the nature of logic and mathematics and of the extensional logic basic to *Principia*, in light of their bearings for epistemology, went along concomitantly" (Ambrose & Lazerowitz 1972: xix). Her system avoids the regress of convention by grounding logic on extension and practice. In other words, Ambrose's dissertation subscribes that extensionalism is the antidote to the epistemic problems within conventionalism and the question of justification in logic, by specifically making the following points:

1. A priori justification is senseless and dependent on human psychology: "Unlike 'asserted' propositions of logic, definitions are neither true nor false; they represent merely an equivalence agreed on… This brings to light the psychological aspect of definition… logic treats definition only as furnishing a symbolic abbreviation for a notion… so far as logic is concerned, definition is verbal" (Ambrose 1932: 18-19).
2. Tautology is logical foundation: "And these laws [identity, contradiction, and excluded middle] are not 'genuine' propositions about material facts. They are tautologies" (Ambrose 1932: 35).
3. Classes are incomplete symbols: "Classes are to be dealt with as symbolic fictions, 'incomplete symbols,' which have no meaning in isolation from the proposition in which they occur" (Ambrose 1932: 48).

4. Formal implication is a subset relation: "The above symbolic equivalence is the sign of the material implication. This, in conjunction with the symbol denoting 'for all values of the variable,' defines formal implication" (Ambrose 1932: 103).

In her 1933 paper "A Controversy in the Logic of Mathematics," which is a substantial revision of the Chapter 4 of her 1932 unpublished PhD dissertation, Ambrose specifically argues that formal systems require interpretation to express anything meaningful. She writes: "According to Hilbert, mathematics has a content independent of logic, a system of extra-logical concrete objects on which intuition falls—but which in themselves are meaningless. These perceptible objects are symbols or 'signs,' but are purged of all intuitive content so as not to be identifiable with the intuitively grasped natural numbers of Brouwer… Mathematics for Hilbert is a language, even something less than a language; each sign, all mathematical axioms and logical principles, the 'elementary deduction' (material implication)—all are rid as far as possible of any meaning. Only the discussion of how this meaningless subject-matter is to be manipulated, i.e., only the 'metamathematics,' which concerns itself with the possibility of certain combinations of signs and the demonstrability of certain theorems about the signs, has meaning" (Ambrose 1933: 474). By making the point that the bare formal symbols carry no semantic content until someone supplies an interpretation, Ambrose thereby contests the epistemic neutrality of formalism, which will not sit well with atomistic postulates such as those defended by Russell and Carnap, for instance.

Briefly, Ambrose's unpublished dissertation and two early published articles repeatedly flag objections to treating logic as merely a matter of social or legislative convention—objections that Margaret Macdonald would elaborate in her 1934 unpublished dissertation and that Quine (1936 and 1951) would later dramatize. Below are the loci in which Ambrose precedes both Macdonald's and Quine's core insights.

### 3. Anticipating Criticism (I): Truth by Convention.

As Ambrose makes clear, it might not be the best procedure to simply declare certain foundational assumptions as true without any proof so a mathematical or logical system can be built. This skepticism is the basis for her criticism of Lewis's intensional system, and the reason for her siding with intuitionism against logicism and formalism when it comes to accepting the requirement of empirical and cognitive grounding of mathematical and logical proof.

Taking as target C. I. Lewis's notion of a priori as a deliberate choice and act of the mind, Ambrose's early publications contend that logical principles cannot be reduced to mere social or legislative conventions—i.e., fiats—because once chosen they become "eternally true" yet evade empirical refutation. Specifically, she writes: "The satisfaction of an expectation by direct verification of meaning, could not make a conceptual system or a priori proposition true. These are eternally true, apart from empirical verification and apart from facts which are instances of them" (Ambrose 1931: 378). And later she will add: "All postulates represent a fiat, discovery then taking place within the given set of postulates. The insistence of some logisticians upon an eternal existence of mathematical facts apart from discovery and upon the eternality of truth about them is very unfortunate… infects the science with subjectivity" (Ambrose 1933: 609). Thus, Ambrose's core thesis—that logic cannot be merely based on a priori concepts and definition—prefigures both Macdonald's and Quine's insistence on empirical requirements and proof in logic.

Ambrose stresses that logical postulates, once chosen, are never strictly falsified by experience—echoing Macdonald's (1934) and Quine's (1936) later point that if they were mere conventions they would require a prior acceptance of themselves. All three authors, thus, deny that logical laws derive legitimacy from social fiat. Specifically speaking about C.I. Lewis's attempt to ground pragmatic truth on the abovementioned conventions, Macdonald (1934) argues that any attempt of founding logical truths on stipulation presupposes the very principles it seeks to establish, leading to some unclear origins. In other words, conventions presuppose logic, not vice versa. MacDonald writes: “It would seem then as though the principles of logic or at least the principle of non-contradiction or self-consistency must be on a different level even from other necessary and a priori propositions… The whole position is very close to that of Kant. We do not agree because we do not have a common world but we have a common world because we agree” (1934: 171, 176).[6] Quine will later maintain that definitions transmit but cannot found truth, since applying definitions presupposes logical inference and, thereby, begs the question of logic’s source. He says: “The difficulty is that if logic is to proceed mediately from conventions, logic is needed for inferring logic from the conventions” (Quine 1936: 258-259).

The convergence between Ambrose, Macdonald, and Quine becomes especially clear once we situate their arguments against the background of early twentieth-century attempts to secure the foundations of logic and mathematics through convention. In the 1920s and 1930s, logical empiricists and pragmatists such as Carnap and C. I. Lewis sought to explain the necessity of logical truths by treating them as products of stipulation—rules we adopt to structure inference, rather than truths discovered or justified independently. This legislative conception of logic promised to reconcile the apparent apriority of logical principles with an empiricist epistemology: Logical truths would be necessary because we make them so. However, the very attractiveness of this view exposed its central vulnerability. For if logical laws derive their authority from conventions, then the act of adopting those conventions itself must presuppose some logical framework. It is precisely this tension that Ambrose identifies earlier and more sharply than has been recognized.

[6] All references to MacDonald (1934), if not indicated, are taken from Spinney (2025).

Ambrose's critique begins from the observation that once logical postulates are adopted, they are treated as unrevisable—"eternally true," as she puts it—yet their supposed conventional origin provides no account of why they should enjoy this privileged status. Her point is not merely epistemic but structural, for a convention cannot ground the very inferential machinery required to formulate, adopt, or apply it. This is why she insists that logical principles cannot be falsified by experience but also cannot be justified by fiat. The act of stipulating a logical rule already presupposes the validity of some inferential norms and, thus, the conventionalist picture collapses into circularity. Ambrose's insight directly anticipates Macdonald's later formulation of the same problem. In her 1934 dissertation, as seen, Macdonald argues that any attempt to found logical truths on stipulation presupposes the principle of non-contradiction or self-consistency, thereby placing logic on a different level from other a priori propositions. Her Kantian-sounding remark—that we have a common world because we agree, not the other way around—captures the same structural point Ambrose had already articulated. That conventions presuppose logic, not vice versa.

Quine's 1936 critique of truth by convention sharpens this regress argument into a general challenge to the conventionalist program. Definitions, he argues, can transmit truth but cannot originate it, because applying a definition requires logical inference. Thus, if logic is to be derived from conventions, logic must already be in place to carry out the derivation. Quine's formulation is more technically explicit, but the core insight is the same one Ambrose had pressed against Lewis. That conventionalists cannot escape presupposing the very logical principles they claim to ground. Taken together, therefore, these three critiques reveal a shared recognition of the structural instability of logical conventionalism at a moment when formalism, pragmatism, and logical empiricism were all attempting to secure the foundations of mathematics through stipulation. Ambrose's contribution is not merely parallel to Macdonald's and Quine's. It is chronologically prior and conceptually aligned, identifying the same regress and the same pre-conditional role of logic before either Macdonald or Quine had written their critiques.

Recall that Ambrose (1931) already argues that stipulative definitions of logical terms inevitably loop back on themselves. She considers that they end up in a vicious regress and circularity. She says: "The view Mr. Lewis is concerned to combat is that there would be an infinite regress if relational terms were defined as relations of relations… This statement hinges on the view that definition is eventually circular, not linear. Our system of concepts as a whole is closed, each meaning to be found in terms of every other. […] I really cannot see any advantage of a circular procedure over a linear regress, either logically or epistemologically" (Ambrose 1931: 376-377). Following this thought, Ambrose underscores that logic cannot be demoted to the status of one more contingent hypothesis. She insists that it must be in place before any mathematical or linguistic conventions can function. In other words, logic is a necessary precondition for the establishment of conventions, circling back to what was established above. She writes: "The acceptance of the law of excluded middle means the belief in the solubility of every problem, since every proposition is viewed as either true or false. The rejection of it probably entails greater precision in the statement of theorems and the confinement of mathematics to a domain of perfect security… The issue is one as to methodology, one which will be settled, not logically, but on pragmatic grounds" (Ambrose 1933: 610).

Ambrose's view of logic's pre-conditional role, therefore, foreshadows both Macdonald's and Quine's arguments. Later, in her dissertation, Macdonald develops the same point regarding circularity. According to her, any legislative act to establish logical truths presupposes those truths, collapsing into an endless regress (MacDonald 1934: 171-174). Quine, as is well known, sharpens this regress objection. Deriving particular truths from general conventions, according to him, requires inference and thus the prior establishment of logical laws, what generates an infinite regress (Quine 1936: 270). Nonetheless, Ambrose insists that without a minimal set of logical principles already in place, no linguistic or mathematical convention could be formulated. In her own words: "If there is to be understanding between two minds, they must have some identical elements in common… Reality is constructed on the basis of concepts 'identified by the relational patterns which speech and behavior in general are capable of conveying.' A and B who understand each other frame concepts in terms of identical systems of orderly relation indicated by behavior, and in consequence find themselves in a common world" (Ambrose 1931: 369). Briefly, unless a minimum set of principles is already operative, no convention in mathematics—or language—could even get off the ground. Macdonald's dissertation concurs with Ambrose that logical principles are prerequisites for any meaningful convention or communication, not their product (1934: 203). While acknowledging that we may stipulate definitions within a chosen framework, Quine later denies that the selection of the framework itself can be conventional, since stipulation already presupposes that framework's inferential rules (1936: 260, 270).

Ambrose specifically highlights what she calls the paradox of the a priori fiat[7]: That logical postulates, besides being the product of "deliberate choices of the mind," yet once adopted they are held as immutably true and forever beyond empirical disproof, making hard to find even their grounds. She says: “A priori truths are never proved false, only inconvenient. And as such they are abandoned. Yet even though abandoned, they are eternally true. Any contradiction ‘new truth’ presents to the old is only verbal, for the old terms have simply acquired new meanings. Concepts may be discarded but cannot change” (Ambrose 1931: 367). Ambrose, therefore, recognizes that what we call logical necessity must necessarily resist apriority, as a mere definitional or conventional grounding, on the faulty consequences of impeding revision of the same principles one established. Echoing this, Macdonald insists that imagining alternative conventions for principles like non-contradiction is absurd, indicating that logical laws occupy a special, necessary status (1934: 171-174). Quine concedes that within a formal system certain truths can be true by stipulation but argues that ultimate logical laws cannot derive their obligatory force from definitions alone (1936: 250), they require some additional empirical discoveries. If fact, as Quine (1951: 40) later emphasizes: “No statement is immune to revision. Revision even of the logical law of the excluded middle has been proposed as a means of simplifying quantum mechanics.”

[7] Spinney (2025: 5) says: "In support of his contention… that Quine's chief target was [C. I.] Lewis, Morris (2017: 371) identifies the continuity in terminology with respect to 'legislate' and its cognates appearing both in *Mind and the World Order* and 'Carnap on Logical Truth' [(Quine 1960b)]. If Morris is correct in this contention, then Macdonald's anticipation of Quine is striking not only in the similarity of argument they employ… but also in their having the same target in mind." If Spinney is correct in his contention that Morris's (2017) terminological continuity, choice of target, and argumentative similarity demonstrates MacDonald's anticipation of Quine's criticism, therefore, I argue, Ambrose anticipates both MacDonald's and Quine's criticism for these same reasons.

In addition, Ambrose’s separation of pragmatic adoption from foundational justification aligns with Quine’s later holistic epistemology. Although she accepts that pragmatic considerations such as simplicity and rigor may influence which logical principles we adopt, she argues that this choice lies outside pure logic itself and cannot render logical truths conventional (Ambrose 1933: 610). Macdonald will similarly later allow that refinement of logical notation may serve pragmatic ends but, she would insist, this does not convert logical origins into conventions (1934: 203). In fact, Macdonald repudiates justification for logical truths. As she later reformulates from her 1934 dissertation: “The suggestion that statements of logical truth admit of justification is senseless” (Macdonald 1940: 26). As is well known, Quine’s later holism incorporates pragmatic and empirical factors in simplicity and decision-making yet cautions that such considerations recommend frameworks without justifying logical truth itself, which remains analytic in name only (1951: 41-43).

**4. Anticipating Criticism (II): The Analytic-Synthetic Dichotomy.**

Ambrose employs the terminology for analyticity and apriority, and their correlates, common during the 1930s indistinctly. Therefore, Ambrose does not posit the analytic-synthetic divide *per se*. Nevertheless, Ambrose's functional distinction between conceptual and factual statements foreshadows Quine's later critique of the analytic-synthetic distinction. Although she does not employ the labels "analytic" and "synthetic" as determined by the dichotomy, Ambrose distinguishes intensional (or conceptual or conventional) truths from extensional (or factual) ones, treating conceptual relations as prior to any factual claims (1931: 373-375). Quine famously pronounces the analytic-synthetic distinction "untenable," as White (1950) refers to it, arguing that no clear, non-circular criterion separates truths of meaning from truths of fact (1951: 20). Both Ambrose and Quine share criticism of formalism in logic and atomism in epistemology, Lewis and Russell being her preferred target (Cf. Ambrose 1931 and 1937)[8] and, if Morris (2017) is right (see note 7), Carnap and Lewis, Quine's.

[8] Russell (1936: 141-145) accuses Ambrose of misunderstanding the analytic-synthetic dichotomy. He specifically says: "Miss Ambrose, if we are to interpret her literally, must hold that the statement 'all men are mortal' is neither true nor false, since we are not 'certain of being able to verify it or prove it false.'…the form of words 'all men are mortal' is outside the scope of the Law of Excluded Middle. […] I hold that, as soon as I know what is meant by 'men' and what by 'mortal,' I know what is meant by 'all men are mortal,' and I know quite certainly that either this statement is true or some man is immortal" (Russell 1936: 145). Chapman (2024) does not mention Russell's response to Ambrose. Loner (2024), even though highlights it, paradoxically only insinuates that Russell is simply ridiculing Ambrose's critiques towards logicism. As mentioned, Russell (1936) clearly indicates that she is directing attacks to the way formalism and logicism understand the aprioristic propositional foundations of logic and mathematics, precisely how it was established by the positivists and that he himself, having been an atomist, still endorses. I, therefore, believe not only that Russell is taking Ambrose seriously, but I have also arrived to believe that he understood as early as in 1936 that Ambrose was already advancing some damaging critics to the sharp dichotomy between the analytic and the synthetic, which Quine will later elevate.

Ambrose specifically refers to the obscurity of the notion of analyticity, at least as defined under pragmatic considerations and, particularly, as sustained by C. I. Lewis. By showing that intensional definitions are vulnerable to the same circularity as extensional ones, Ambrose already exposes the fragility of any "analytic" category grounded in meaning alone. Beyond that, Ambrose writes: "That concepts are got by abstraction and then held either as tentative or definitive is an account upon which Mr. Lewis makes his position especially clear. […] In one surprising passage, Mr. Lewis supports a view… obliterating the difference between a priori propositions and empirical generalizations… A stone, for example, would be defined as a freely falling body for certain purposes, while for others such a statement would be merely an empirical generalization. The question which comes first… now becomes insignificant" (Ambrose 1931: 373-374). Ambrose, once again, is underscoring that logical laws, once adopted, carry the force of necessity even though they could have begun as mental fiat (Cf. Ambrose 1931: 367; Ambrose 1933: 609, above). Quine (1951) will later argue that attempts to define analyticity via synonymy or interchangeability presuppose analytic truths, thus failing to demarcate the boundary non-circularly.

Famously, he writes: "Our argument is not flatly circular, but something like it. It has the form, figuratively speaking, of a closed curve in space. Interchangeability *salva veritate*, as we formulated it, presupposes that the language admits an intensional 'necessarily' whose job is quietly to flag analytic truths. But can we condone a language which contains such an adverb? Does the adverb really make sense? To suppose that it does is to suppose that we have already made satisfactory sense of 'analytic.' Then what are we so hard at work on right now?" (Quine 1951: 29). Quine, thus, has no doubt that any attempt to define analyticity ends up assuming the very notion it is trying to explain. His test for analyticity—substitution *salva veritate* using a modal operator like "necessarily"—already presupposes that we understand which truths are analytic. The argument is, therefore, not a simple circle but a "closed curve," for every route we take to define analyticity loops back to analyticity itself. If we allow "necessarily" into the language then, we are already assuming the distinction we are trying to justify, which makes the entire project self-defeating. Ambrose's (1931 and 1933) concerns about circular definitions, hence, anticipate Quine's more famous argument against the analytic-synthetic dichotomy.

For, even though Ambrose's work prefigures Quine's holistic overturning of the analytic-synthetic divide, she employs the language of extensional versus intensional or formal logic and will later link this to the normativity embedded in ordinary language, like also would MacDonald. By emphasizing logic's pre-conditional status though, Ambrose implicitly rejects an atomistic view where individual "analytic" truths stand apart from the rest of knowledge. Like Macdonald and Quine would do later, Ambrose's appeal to a regress objection entails, as seen above, that logical truth cannot be isolated for convention without invoking the whole system. Quine makes this explicit in "Two Dogmas," arguing for confirmation holism. He specifically says: "The totality of our so-called knowledge or beliefs, from the most casual matters of geography and history to the profoundest laws of atomic physics or even of pure mathematics and logic, is a man-made fabric which impinges on experience only along the edges. […] No statement is immune to revision. Revision even of the logical law of the excluded middle has been proposed as a means of simplifying quantum mechanics. […] Our statements… face the tribunal of sense experience only as a corporate body" (1951: 39-41).

Ambrose's restrained view of pragmatic revision, therefore, marks an intermediate stance. She upholds a meaningful conceptual distinction even while recognizing logical norms as human choices. She, thus, allows that pragmatic needs may lead us to adopt or abandon certain logical norms, but insists that such revision does not blur the conceptual distinction between conceptual truths and empirical claims (Ambrose 1933: 610). Quine, in contrast, pretends to dissolve the analytic-synthetic distinction here altogether, besides recognizing analyticity later, by holding that revisions in any part of our web of beliefs may shift our acceptance of formerly analytic truths (Quine 1951: 41). It is not then surprising that Quine also adopted extensionalism as a philosophical orientation, allowing him to put further distance from intensional logic and stipulative definitions (Cf. Quine 1960a and, especially, Quine 2001).

Ambrose (1935a and 1935b) further develop her skepticism toward non-verifiable logical inference. Even though the aim of these papers is to defend finitism against objections, she calls for epistemic responsibility in mathematical reasoning. She argues: “The finitist demands that we should be certain of being able to verify or to prove false a verbal form before we hold it to be either true or false in any clear sense of these two words” (Ambrose 1935a: 379). Briefly, Ambrose here claims moderate finitism as epistemic criticism. But this is a matter for another occasion.

**5. The Macdonald-Ambrose Partnership at Cambridge (1934-1935).**

Ambrose arrived in Cambridge with a Wellesley post-doctoral fellowship in summer of 1932, which soon proved insufficient.[9] Given the impossibility for women to occupy certain academic positions, even though she had already completed a PhD from the University of Wisconsin-Madison, she was obligated to enroll as a PhD student to attend the Moral Sciences meetings and other events in exchange for some economic aid. She was part of Wittgenstein’s inner circle of students,[10] attending his private lectures and dictations until her departure before the beginning of the 1935 fall term to occupy a position as lecturer in Mathematics at the University of Michigan. Ambrose transcribed many of Wittgenstein's remarks in detail, work that Macdonald consulted and contributed to before Wittgenstein objected to their publication (Wittgenstein 1979).

[9] It is worth noting that this is the same fellowship Quine received that same year, which allowed him to visit Vienna to study under Carnap and, the summer after, have a short stay with Tarski in Warsaw. It is not unlikely that, been part of the same cohort, Quine and Ambrose met at some point prior of their respective departure. In an anecdotical note, Ambrose declined Oberlin College’s offer to their graduate Philosophy program in 1928, where Quine was a Philosophy undergraduate, taking Wisconsin’s offer instead (Millikin University’s Special Collections, Alice Loman Ambrose Archive, Box 1, folder 1, folio 16, CV). I am not aware of any information about her visiting the college, but if that happened, it is likely that Quine and Ambrose met then.

Macdonald's 1934 unpublished doctoral thesis, supervised by Susan Stebbing and defended at University College London, addresses the question of logical necessity and truth, but centers on the requirements for linguistic meaning and communication, themes deeply consonant with Ambrose's concerns. She arrived in Cambridge during 1934, although it is likely that she had attended the Moral Sciences meetings prior, where she likely met Ambrose. She did not obtain permission to participate in Wittgenstein's classes nor dictations, and there is no indication that Ambrose read or helped MacDonald with her dissertation. Their friendship and shared methodological background, nevertheless, influenced their mutual skepticism of rigid foundationalism. Macdonald's later work on necessary propositions, framed in Wittgensteinian language philosophy, for example, mirrors Ambrose's constructive critique of logical form and meaning.

There is no doubt that Alice Ambrose and Margaret Macdonald were close intellectual collaborators during their overlapping time at Cambridge in 1934 and 1935. Nevertheless, one very important thing to remember and that many seem to forget, including those that are highly sympathetic to Ambrose's philosophical contributions, is that Ambrose arrives in Cambridge after already completing and defending a PhD dissertation at the University of Wisconsin-Madison.

---

[10] She is Wittgenstein only known student, at least for the roughly couple of years that she was enrolled under his supervision, until, the word is, she dear to publish some of Wittgenstein's early view on the foundations of mathematics, which Wittgenstein considered "quite indecent," as he wrote in a letter to Moore (Wittgenstein 1936). The work in question is Ambrose (1935a and 1935b), which Ambrose decided to publish, without Wittgenstein's final approval, encouraged by Moore, who had taken over the role of editor at *Mind*.

This early work, as described above, is focused on philosophical logic and centers on a defense of the benefits of an extensional logic over intensional logical systems. It is only because women could not occupy certain academic roles, likely due to the well-documented misogynist climate of the era (Cf. Connell & Janssen-Lauret 2022; Janssen-Lauret 2024), that she is obligated to pursue a second doctoral degree, which she will complete and defend in 1938 under the final supervision and approval of G.E. Moore. Although both women would diverge in style, Ambrose tending toward formal analysis while Macdonald diverging more toward ordinary language, their shared Cambridge education and philosophical sensibilities created a fertile environment for the evolution of the foundational critique, one that with no doubt Quine witnessed.

**7. Conclusion**

Ambrose (1931) and (1933) identify and name the fatal weaknesses of any view that would legislate logic into being: The a priori fiat, circularity/regress, and the pre-conditional status of logical laws. Macdonald's 1934 dissertation picks up each thread—denying that logical truths are conventions, exposing the infinite regress in attempts to ground them on agreement, and insisting on logic's priority to convention—thus anticipating Quine's (1936) anti-conventionalism.

Alice Ambrose, therefore, deserves recognition not merely as a transcriber of Wittgenstein or critic of C. I. Lewis, but as a foundational thinker whose work straddles logic and epistemology and criticizes conventionalism. Her publications from the early 1930s prefigure many of the themes that Macdonald's 1934 dissertation would refine and, later, Quine would canonize. The affinities across their separate criticisms suggest not three isolated voices, but a shared intellectual atmosphere that began dismantling the notion that logic's truth could be justified by convention.[11]

[11] Acknowledgments.